\documentclass[a4paper]{amsart}
\usepackage{amsmath}
\usepackage{amssymb}
\usepackage{hyperref}
\usepackage{xcolor}
\usepackage{graphicx}
\begin{document}
\title[Elastic shape matching]{Inexact elastic shape matching in the square root normal field framework}
\thanks{P.~Harms is supported by the Freiburg Institute of Advanced Studies in the form of a Junior Fellowship. N.~Charon is supported by the NSF grant 1819131.}
\author{Martin Bauer}
\address{Martin Bauer: Florida State University, Tallahassee FL 32304, USA}
\email{bauer@math.fsu.edu}
\author{Nicolas Charon}
\address{Nicolas Charon: John Hopkins University, Baltimore MD 21218, USA}
\email{charon@cis.jhu.edu}
\author{Philipp Harms}
\address{Philipp Harms: Albert-Ludwig-University Freiburg, 79104 Freiburg, Germany}
\email{philipp.harms@stochastik.uni-freiburg.de}
\maketitle 
\begin{abstract}
This paper puts forth a new formulation and algorithm for the elastic matching problem on unparametrized curves and surfaces. Our approach combines the frameworks of square root normal fields and varifold fidelity metrics into a novel framework, which has several potential advantages over previous works. First, our variational formulation allows us to minimize over reparametrizations without discretizing the reparametrization group. Second, the objective function and gradient are easy to implement and efficient to evaluate numerically. Third, the initial and target surface may have different samplings and even different topologies. Fourth, texture can be incorporated as additional information in the matching term similarly to the {\tt fshape} framework. We demonstrate the usefulness of this approach with several numerical examples of curves and surfaces. 
\end{abstract}
\section{Introduction}

\subsection{Context.}
The statistical analysis of datasets of curves and surfaces is an active research field with many applications in e.g.\@ computer vision, robotics, and medical imaging; see \cite{younes2010shapes,bauer2014overview,srivastava2016functional} and references therein. 
A recurring and fundamental task is finding optimal point correspondences between given shapes (i.e., the matching or registration problem), where optimality is typically expressed in terms of an elastic deformation energy.   
Solving the elastic matching problem in a numerically efficient way, which scales well to high-dimensional data encountered in real-world applications, remains a major challenge to date.  

\subsection{Relation to previous work.}
This paper draws on two lines of work: 
square root normal fields (SRNFs) \cite{srivastava2011shape,kurtek2012elastic,jermyn2012elastic,jermyn2017elastic}, which allow one to efficiently calculate elastic distances between parametrized shapes, 
and varifold distances \cite{vaillant2005surface,glaunes2008large,charon2013varifold,roussillon2016kernel,kaltenmark2017general}, which are distances between unparametrized shapes without any elastic interpretation.
For each of these frameworks, efficient numerical implementations have been developed.  

\subsection{Contribution.}
We propose a new algorithm which combines SRNFs with varifold distances and inherits many advantages of both approaches. 
The key idea is to use varifold distances to relax the terminal constraint in the elastic matching problem. 
This bypasses the discretization of the reparametrization group, thereby eliminating the main computational burden in previous implementations of SRNF-based elastic shape matching. 
The resulting optimization problem is easy to implement and yields good results on some preliminary experiments on curves and surfaces. 
Moreover, the varifold distances allow one to match shapes with different meshes and even different topologies and to use texture information as in the {\tt fshape} framework. 

\section{Shape analysis of curves and surfaces}

\subsection{Elastic shape analysis.}

Elastic shape analysis operates in a Riemannian framework where infinitesimal shape deformations are measured by a Riemannian metric, which is often related to an elastic (or plastic) deformation energy; see the surveys \cite{bauer2014overview,jermyn2017elastic}. 
We consider parameterized shapes as elements of the Fr\'echet manifold $\operatorname{Imm}(M,\mathbb R^d)$ of immersed hypersurfaces of a $(d-1)$-dimensional compact manifold $M$ into $\mathbb R^d$. 
The corresponding space of unparameterized shapes is the quotient space $B_i(M,\mathbb R^d)=\operatorname{Imm}(M,\mathbb R^d)/\operatorname{Diff}(M)$, 
whose elements are denoted by $[f]=\{f\circ\varphi; \varphi \in \operatorname{Diff}(M)\}$.
Given a $\operatorname{Diff}(M)$-invariant weak Riemannian metric $G$ on $\operatorname{Imm}(M,\mathbb R^d)$, one defines a pseudo-distance between any two immersions $f_0,f_1 \in \operatorname{Imm}(M,\mathbb R^d)$ and their equivalence classes $[f_0],[f_1]\in B_i(M,\mathbb R^d)$ by
\begin{gather}
\label{equ:dist_imm}
\operatorname{dist}_{\operatorname{Imm}}(f_0,f_1)^2
=
\inf_{\substack{f\in C^\infty([0,1],\operatorname{Imm}(M,\mathbb R^d))\\f(0)=f_0, f(1)=f_1}}\ \int_0^1 G_f(\partial_t f,\partial_t f) dt,
\\
\label{equ:dist_bi}
\operatorname{dist}_{B_i}([f_0],[f_1])
=
\inf_{\varphi\in \operatorname{Diff}(M)} \operatorname{dist}_{\operatorname{Imm}}(f_0,f_1\circ\varphi).
\end{gather}
Symmetry of the pseudo-distance on $B_i(M,\mathbb R^d)$ follows from the invariance of the metric $G$ with respect to reparametrizations.
Under suitable conditions on $G$ the pseudo-distance is a distance, i.e., it separates points in the shape space of curves \cite{michor2007overview,mennucci2008properties} or surfaces \cite{bauer2011sobolev}.

From a numerical perspective, the challenge is to calculate the above distances and the corresponding optimizers efficiently. 
This minimization can be solved numerically by path straightening and geodesic shooting methods (see e.g.~\cite{huang2016riemannian,bauer2017numerical}) or as in the next section by exploiting isometries to simpler spaces.

\subsection{Square root normal fields.}

Problem \eqref{equ:dist_imm} simplifies considerably for certain first order Sobolev metrics \cite{younes1998computable,younes2008metric,srivastava2011shape,sundaramoorthi2011new,kurtek2012elastic,jermyn2012elastic,bauer2014constructing}. 
One class of such metrics is defined using square root normal fields (SRNFs), which were introduced by Srivastava e.a.\@ \cite{srivastava2011shape,kurtek2012elastic} for planar curves and later generalized to surfaces by Jermyn e.a.\@ \cite{jermyn2012elastic}. 
The SRNF of an oriented immersed hypersurface $f \in \operatorname{Imm}(M,\mathbb R^d)$ is defined as $\smash{\tilde n_f = n_f \operatorname{vol}_f^{1/2}}$, where $n_f$ is the unit normal field and $\smash{\operatorname{vol}_f^{1/2}}$ the half density of $f$.
For example, the SRNF of a planar curve $f\in\operatorname{Imm}(S^1, \mathbb R^2)$ is given in coordinates $\theta \in S^1$ by $\smash{\tilde n_f=if_\theta \|f_\theta\|_{\mathbb R^2}^{-1/2}}$, where $i$ denotes rotation by 90 degrees and coordinates in subscripts denote derivatives. 
Similarly, the SRNF of a surface $f\in\operatorname{Imm}(S^2, \mathbb R^3)$ is given in coordinates $(u,v) \in S^2$ as $\smash{\tilde n_f=(f_u\times f_v) \|f_u\times f_v\|_{\mathbb R^3}^{-1/2}}$. 
In general, one obtains an elastic pseudo-Riemannian metric on $\operatorname{Imm}(M,\mathbb R^d)$ by setting
\begin{equation*}
G_f(h,k) 
=
\int_M \langle D_{(f,h)}\tilde n_f,D_{(f,k)}\tilde n_f\rangle_{\mathbb R^d},
\end{equation*}
where $D_{(f,h)}\tilde n_f$ denotes the directional derivative of $\tilde n_f$ at $f$ in the direction $h$.
This pseudo-Riemannian metric $G$ is $\operatorname{Diff}(M)$-invariant, and by construction the map $f\mapsto\tilde n_f$ is a Riemannian isometry into the flat space of square integrable vector-valued half densities.
For curves one obtains a Riemannian metric by modding out translations.
For surfaces the situation is more complicated, as described in \cite{jermyn2012elastic}, and the kernel of the pseudo-metric may be larger than only translations.
The metric belongs to the class of first order Sobolev metrics, which have been studied in great detail \cite{michor2007overview,mennucci2008properties,bauer2011sobolev}.

The advantage of this construction is that the Riemannian distance of $G$ on $\operatorname{Imm}(M,\mathbb R^d)$ can be approximated efficiently as follows:
\begin{equation}\label{equ:approx}
\operatorname{dist}_{\operatorname{Imm}}(f_0,f_1)
\approx
\|\tilde n_{f_0}-\tilde n_{f_1}\|_{L^2}.
\end{equation} 
Equality holds whenever the straight line between $\tilde n_{f_0}$ and $\tilde n_{f_1}$ is contained in the range of the map $f\mapsto\tilde n_f$.
In general, equality holds up to first order for $f_0$ close to $f_1$ because the map $f\mapsto\tilde n_f$ is a Riemannian isometry. 

The approximate distance \eqref{equ:approx} descends to the quotient space $B_i(M,\mathbb R^3)$ as described in \eqref{equ:dist_bi}.
However, \eqref{equ:dist_bi} involves a minimization over the reparametrization group, which is computationally costly.
For curves this can be solved by dynamic programming \cite{srivastava2011shape} or using an explicit formula \cite{lahiri2015precise}. 
For surfaces in spherical coordinates, Jermyn~e.a.~ \cite{jermyn2012elastic} proposed to discretize the diffeomorphism group of the two-dimensional sphere using spherical harmonics. 
This article puts forth an alternative method for minimization over the reparametrization group, which is based on varifold distances.

\subsection{Varifold distances.}

Geometric measure theory provides several embeddings of shape spaces into Banach spaces of distributions \cite{vaillant2005surface,glaunes2008large,charon2013varifold,roussillon2016kernel,kaltenmark2017general} with corresponding metrics. Varifold embeddings are one instance of this construction and are defined as follows (cf. \cite{kaltenmark2017general} for details).
Given a reproducing kernel Hilbert space $W$ of real-valued functions on $\mathbb R^d\times S^{d-1}$, 
one associates to any immersion $f \in \operatorname{Imm}(M,\mathbb R^d)$ the varifold $\mu_f\in W^*$ which satisfies
\begin{equation*}
\forall w \in W: 
\qquad 
(\mu_f | w )_{W^*,W}
=
\int_M w(f(x),n(x)) \operatorname{vol}_f(dx).
\end{equation*}
The map $f\mapsto \mu_f$ is reparametrization-invariant and, under suitable assumptions on the kernel of $W$, injective \cite{kaltenmark2017general}. 
Thus, one obtains a well-defined distance on the quotient space $B_i(M,\mathbb R^d)$ by defining for any two immersions $f_0,f_1 \in \operatorname{Imm}(M,\mathbb R^d)$:
\begin{equation*}
\operatorname{dist}_{\operatorname{Var}}([f_0],[f_1])
=
\|\mu_{f_0}-\mu_{f_1}\|_{W^*}.
\end{equation*}
From a computational point of view, these distances have explicit expressions in terms of the kernel function of $W$ and are easy to implement for discrete curves or surfaces. 
We will use such distances to relax the terminal constraint in the boundary value problem for geodesics on shape space, as described next.

\subsection{Combining SRNFs and varifold distances.}

Square root normal fields and varifold distances can be combined in an efficient matching algorithm for unparametrized shapes.
This idea has been previously used in combination with large deformation models in \cite{charon2013varifold,kaltenmark2017general} and with $H^2$ metrics on the space of curves in \cite{bauer2018relaxed}.
The boundary value problem \eqref{equ:dist_bi} for geodesics on $B_i(M,\mathbb R^d)$ can be formulated as the program
\begin{equation}
\underset{f}{\text{minimize}}
\hspace{1ex} 
\operatorname{dist}(f_0,f)
\qquad
\text{subject to}
\hspace{1ex}
\operatorname{dist}_{\operatorname{Var}}([f],[f_1])=0.
\end{equation}
Relaxation using a (large) Lagrange multiplier $\lambda$ and approximation of the elastic distance as in \eqref{equ:approx} yields
\begin{equation}\label{equ:our_algorithm}
\underset{f}{\text{minimize}}
\hspace{1ex} 
\|\tilde n_{f_0}-\tilde n_f\|_{L^2}^2 + \lambda \operatorname{dist}_{\operatorname{Var}}([f],[f_1])^2.
\end{equation}
This program has several advantages over previous alternative formulations of the SRNF matching problem \cite{srivastava2011shape,kurtek2012elastic,jermyn2012elastic}. 
First, the objective function and its gradient are easy to implement and can be computed efficiently. 
Second, the initial and target surface may have different discretizations and even different topologies.
Third, texture information can be incorporated into the varifold matching term similarly to the {\tt fshape} framework \cite{charon2014functional,charlier2017fshape}. 

\begin{figure}[t]
\centering
\includegraphics[width=.48\textwidth]{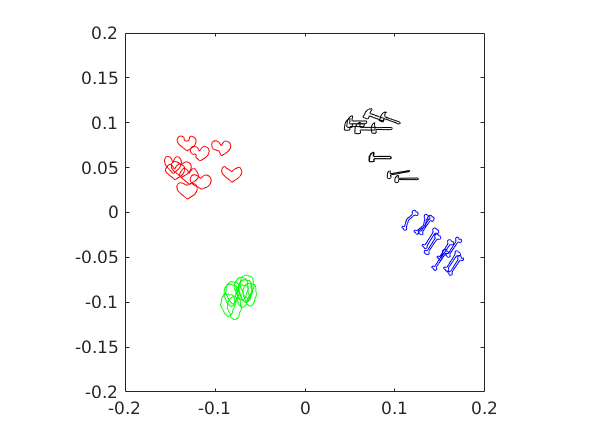}
\includegraphics[width=.48\textwidth]{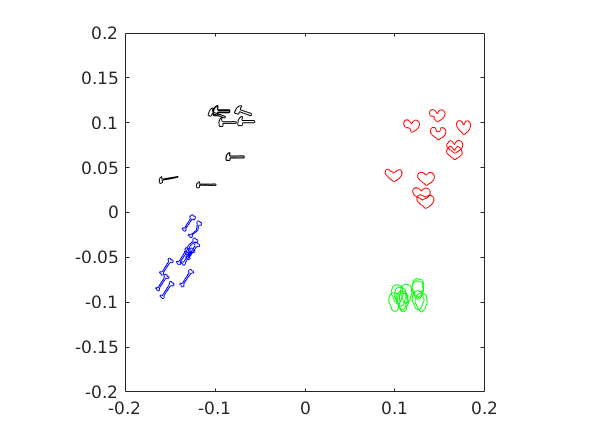}\\
\includegraphics[width=.48\textwidth]{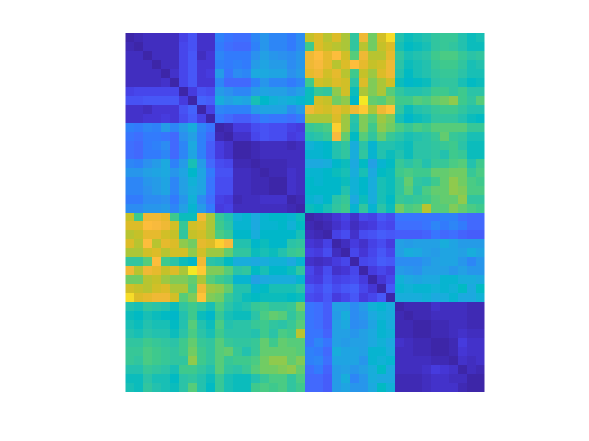}
\includegraphics[width=.48\textwidth]{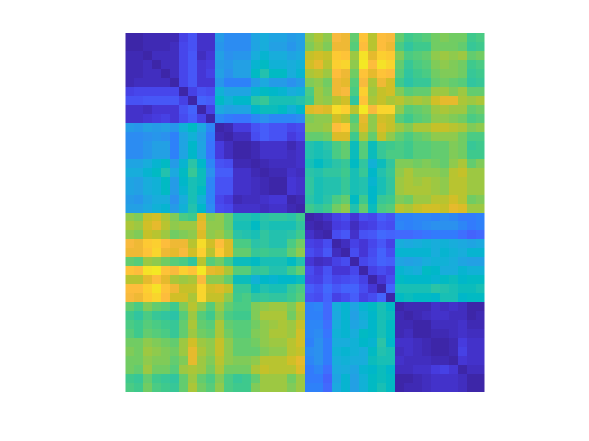}
\caption{Distances and clusters produced by our algorithm are comparable to state-of-the-art curve matching using dynamic programming \cite{srivastava2011shape} when tested on curves in the Kimia database. Left: our SRNF-varifold algorithm; right: dynamic programming; top: distance-based multi-dimensional scaling; bottom: symmetrized distance matrix.}
\label{fig:comparison}
\end{figure}
\section{Numerical implementation and results}

\subsection{Algorithm.}
Given a pair $(f_0,f_1)$ of curves or surfaces, the program \eqref{equ:our_algorithm} looks for a minimizer $f$ with $\operatorname{dist}_{B_i}([f_0],[f_1])=\operatorname{dist}_{\operatorname{Imm}}(f_0,f)$ and $[f]=[f_1]$.
Thus, the algorithm solves the registration problem and calculates the distance between the unparametrized shapes $[f_0]$ and $[f_1]$. 
Note that it does, however, not provide a geodesic homotopy between these shapes. 
Such a homotopy can be obtained from the linear homotopy between $\tilde n_{f_0}$ and $\tilde n_{f_1}$ by (approximate) inversion of the SRNF map $f\mapsto \tilde n_f$.
For open curves this inversion is exact and easy to implement.
For closed curves, the range of the SRNF map is not convex, and an approximate inverse has to be used \cite{srivastava2011shape}. 
For surfaces, this is a delicate issue \cite{jermyn2012elastic}, and to the best of our knowledge there exists no publicly available implementation for general triangulated surfaces.

\subsection{Implementation.}
To implement the program \eqref{equ:our_algorithm} numerically, one has to discretize the space of parametrized shapes. 
An advantage over \cite{jermyn2012elastic} is that the reparametrization group does not need to be discretized.
Piecewise linear curves and triangular meshes are suitable discretizations in our context, the reason being that square root normal fields and kernel-based varifold distances extend naturally to these spaces.
The minimization is performed using an L-BFGS method. 
The gradient of the discretized energy functional \eqref{equ:our_algorithm}, which is needed by the L-BFGS method, has an explicit form and can be implemented efficiently.

\subsection{Curves.}
\begin{figure}[h]
\centering
\includegraphics[trim={4cm 0 4cm 0},clip,width=\textwidth]{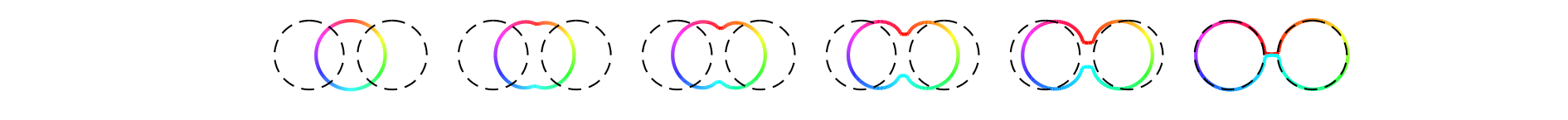}
\\
\includegraphics[trim={4cm 0 4cm 0},clip,width=\textwidth]{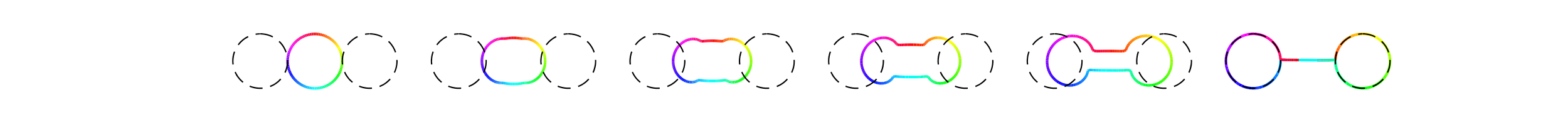}
\caption{Our algorithm can match shapes with different topologies. Left to right: geodesic interpolation (colored) between a single circle and a pair of circles (black, dashed); top: small distance between the pair of circles; bottom: small distance between the pair of circles. }
\label{fig:circles}
\end{figure}
\begin{figure}[h]
\centering
\includegraphics[trim={2cm 0 2cm 0},clip,width=.32\textwidth]{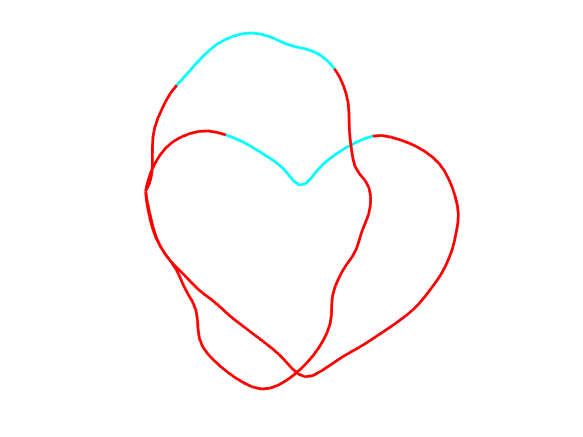}
\includegraphics[trim={2cm 0 2cm 0},clip,width=.32\textwidth]{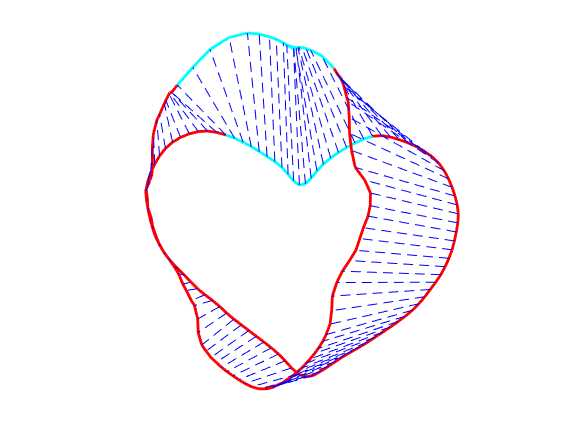}
\includegraphics[trim={2cm 0 2cm 0},clip,width=.32\textwidth]{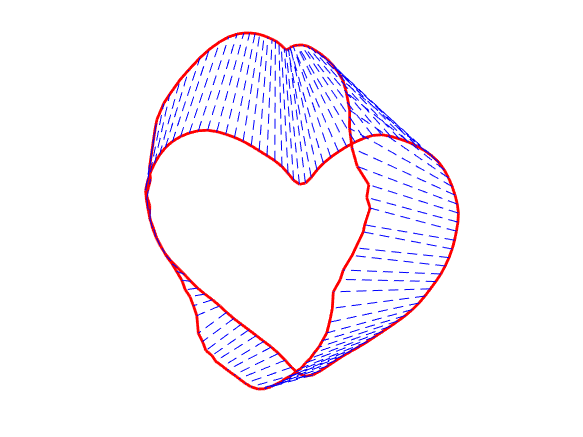}\\
\caption{Elastic matching of curves with functional data. Left: source and target curves with binary functional data in red/cyan. Middle: matching using functional data. Right: purely geometrical matching without functional data.}
\label{fig:fshape_curve}
\end{figure}
For curves, our algorithm is comparable to state of the art methods. 
On the Kimia dataset\footnote{Computer Vision Group at LEMS at Brown University: Database of 99 binary shapes. \url{https://vision.lems.brown.edu/content/available-software- and-databases}} it produces distances and clusters which are similar to those based on dynamic programming, as shown in Fig~\ref{fig:comparison}. 
A nice feature of our algorithm, which stems from the use of varifold distances, is that the initial and target shapes are allowed to have different topologies. 
For example, one can match a single circle to a pair of circles, as demonstrated in Fig.~\ref{fig:circles}.
This is not possible using previous methods for shape matching using SRNFs or SRVTs. 
There are potential applications in cell division and removal of topological noise. 
Another feature of our algorithm is that it can account for functional data on the given shapes, as demonstrated in Fig.~\ref{fig:fshape_curve}.
To this aim, the varifold distance in \eqref{equ:our_algorithm} is replaced by a functional shape distance, as developed in   \cite{charon2014functional,charlier2017fshape}. 
This has several applications.
The functional data may be dictated by the application at hand, as e.g.\@ in the case of texture information. 
An interesting alternative to be explored in future work is to use shape descriptors as functional data to guide the matching algorithm. 

\subsection{Surfaces.}
\begin{figure}
\centering
\includegraphics[trim={2cm 0 2cm 0},clip,width=.48\textwidth]{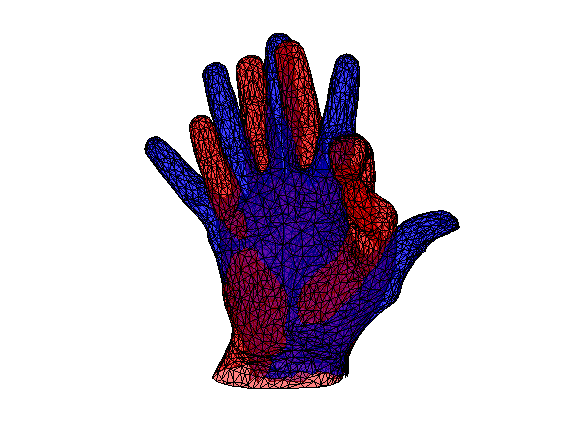}
\includegraphics[trim={2cm 0 2cm 0},clip,width=.48\textwidth]{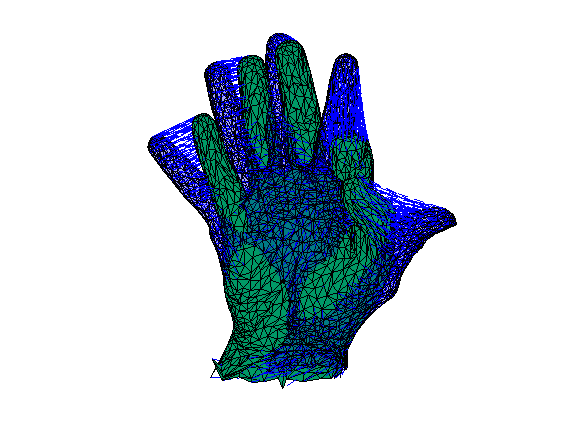}\\
\caption{Anatomically correct correspondences obtained by elastic matching of two surfaces. Left: template $f_0$ (blue, 2322 vertices) and target $f_1$ (red, 2829 vertices). Right: output $f$ of the matching algorithm (green) and a linear homotopy between $f_0$ and $f$.}
\label{fig:hands}
\end{figure}

For surfaces, we obtain some promising first results and see a high potential of improvement over alternative methods. 
An example is presented in Fig.~\ref{fig:hands}, where the optimal point correspondences between two hand postures were calculated. 
As the two triangulated surfaces in this experiment had different mesh connectivities, and no point-to-point correspondences were initially available, we had to initialize the optimization procedure with the template surface.
After optimization using an adaptive choice of Lagrange multiplier $\lambda$ in \eqref{equ:our_algorithm}, we obtained an excellent fit of the deformed template onto the target with anatomically correct point correspondences. 

\bibliographystyle{abbrv}

\end{document}